\documentclass[11pt]{amsart}
\usepackage[latin1]{inputenc}
\usepackage{amsmath,amsfonts,amssymb}
\usepackage{amscd,graphicx}
\usepackage[all,line]{xy}
\usepackage{mathrsfs}
\newcommand{\la}{\lambda}

\newcommand{\N}{{\mathbb{N}}}
\newcommand{\R}{{\mathbb{R}}}

\newcommand{\beq}{\begin{equation*}}
\newcommand{\eeq}{\end{equation*}}
\newcommand{\be}{\begin{equation}}
\newcommand{\ee}{\end{equation}}
\newcommand{\bdm}{\begin{displaymath}}
\newcommand{\edm}{\end{displaymath}}

\newtheorem{theorem}{Theorem}[section]

\newtheorem{prop}[theorem]{Proposition}
\newtheorem{cor}[theorem]{Corollary}
\newtheorem{lemma}[theorem]{Lemma}
\newtheorem{example}[theorem]{Example}
\newtheorem{rem}[theorem]{Remark}

\newtheorem{definition}[theorem]{Definition}

\long\def\symbolfootnote[#1]#2{\begingroup\def\thefootnote{\fnsymbol{footnote}}
\footnote[#1]{#2}\endgroup}

\begin{document}

%\pagestyle{myheadings}

%\markright{\today}
\title{Generalized hook lengths in symbols and partitions}

%\iffalse{===
\author{Christine Bessenrodt}

\address{Institut f\"ur Algebra, Zahlentheorie und Diskrete Mathematik
\\
Leibniz Universit\"at Hannover, Welfengarten 1
\\
D-30167 Hannover, Germany }
\email{bessen@math.uni-hannover.de}

\author{Jean-Baptiste Gramain}
\address{Institut de Mathématiques de Jussieu,
Université Denis Diderot, Paris VII \\
UFR De Mathématiques, 2 Place Jussieu, F-75251 Paris Cedex 05, France
}
\email{gramain@math.jussieu.fr}

\author{J\o rn B. Olsson}
\address{
Department of Mathematical Sciences, University of Copenhagen\\
Universitetsparken 5,DK-2100 Copenhagen \O, Denmark
}
\email{olsson@math.ku.dk}

\date{January 25, 2011}

%===}\fi
\maketitle

\begin{abstract}
In this paper, we present, for any integer $d$, a description of the set of hooks in a $d$-symbol. We then introduce generalized hook length functions for a $d$-symbol, and prove a general result about them, involving the core and quotient of the symbol. We list some applications, for example to the well-known hook lengths in integer partitions. This leads in particular to a generalization of a relative hook formula for the degree of characters of the symmetric group discovered by G. Malle and G. Navarro in \cite{MaNa}.

\end{abstract}

\symbolfootnote[0]{2000 Mathematics Subject Classification 20C30,
  20C33, 05A17}

\symbolfootnote[0]{Keywords: Symbols, Hooks, Hook Lengths, Partitions, Core, Quotient}

%\bigskip
%\begin{center}
%\large{\bf{Introduction}}
%\end{center}

In his work \cite{Ma} on unipotent degrees in reflection groups
G.\ Malle used $d$-symbols as labels and defined hooks in $d$-symbols
and associated a length to a hook. With these he was able to prove a
``hook formula'' for the degrees.

In this paper we introduce generalized hook length functions for
$d$-symbols and prove a general result about them.
More specifically we consider certain functions $h$
 from the set $H(S)$ of hooks of a $d$-symbol $S$
to $\mathbb{R}$ and decompositions of the multiset $\mathcal{H}(S)$ of all
generalized hook lengths $h(z)$,
where $z \in H(S)$. For a given positive integer $\ell$ we
find a decomposition which is compatible with the
$\ell$-core and $\ell$-quotient of $S$.

Based on a crucial well controlled correspondence between hooks
in symbols and hooks in associated core and quotient symbols,
our main result (Theorem~\ref{key})
on generalized hook lengths for symbols is deduced.
This may seem quite special, but in view of its consequence, Theorem~\ref{metasatz},  
it has a number of applications.
For instance, we show that the relative hook formula obtained by
Malle and Navarro  \cite[Theorem 9.1]{MaNa} is just
the well-known hook formula for the degree of the irreducible
characters of the symmetric groups
with the hooks suitably arranged (Remark \ref{managen}).
If $\mathcal{H}(\la)$ is the multiset of hook lengths for a partition
$\la$, and $\la$ has $d$-core partition $\la_{(d)}$ then we have in particular
$\mathcal{H}(\la_{(d)}) \subset \mathcal{H}(\la)$. Furthermore
the remaining elements of $\mathcal{H}(\la)$
may be seen as {\it modified hook lengths} of a $d$-quotient partition for $\la$  
(Theorem \ref{partmeta2}).

\smallskip

Our paper is organized as follows. In section~\ref{betad} we explain
the relations between partitions, $\beta$-sets and $d$-symbols, their
hooks and their corresponding cores and quotients. In section~\ref{hookscor}
we set up useful bijections between the set of hooks of a symbol and the
set of hooks of its core and quotient.  This is used in section~\ref{hookfunc}
to prove our key result, Theorem~\ref{metasatz}. The
theorem is then applied to partitions in section~\ref{appart} and to
symbols in section~\ref{symbs}.
\medskip

\section{On $\beta$-sets and $d$-symbols}\label{betad}

%\medskip

A $\beta$-set $X$ is a  finite subset of $\N_0$.
For $s \in \N_0$ we put $X^{+s}=(X+s)\cup \{0,1,\ldots,s-1\}$.
Write $X=\{a_1,a_2,\ldots,a_t\}$ where we always assume $a_1>a_2>\cdots>a_t$.
Then we associate to the $\beta$-set $X$
the partition $p(X)$ having the non-zero numbers among $a_i-(t-i)$,
$i=1, \ldots,t$, as parts.
For example, $p(\{5,3,0\})=(3,2)$.
Note that $p(X)=p(X^{+s})$ for all $s \in \N_0$.

\smallskip

Let $d \in \N$. We define a $d$-{\it symbol}
$S=(X_0,X_1,\ldots,X_{d-1})$
as a $d$-tuple of $\beta$-sets. In analogy with
$\beta$-sets we define
$S^{+s}=(X_0^{+s},X_1^{+s},\ldots,X_{d-1}^{+s})$.

\smallskip

In \cite{Ma}, a $d$-symbol is an {\em equivalence class}
of such $d$-tuples. The equivalence
relation is generated by the following operations:

$\bullet$  cyclic permutation of the $\beta$-sets

$\bullet$  replacing $S$ by $S^{+1}$.

The multiset of all hook lengths in $S$, as defined in \cite{Ma},
only depends on the equivalence class of $S$. The more general
definition of hook lengths in symbols which we introduce in section 3
includes those of \cite{Ma}, but they are
not independent of the cyclic permutations of the $\beta$-sets.

\smallskip

Let $s \in \N$.  We set $[s]=\{0,1,\ldots,s-1\}$.
For $n \in \mathbb{N}_0$, we define $n_{[s]} \in [s]$ as the
remainder of $n$ after division by~$s$.

\begin{definition}\label{asssymb}
Given  $d \in \N$ we associate to a $\beta$-set $X$ a $d$-symbol
$$s_d(X)=(X^{(d)}_0,X^{(d)}_1,\ldots,X^{(d)}_{d-1})$$
where for $j \in [d]$
$$X^{(d)}_j=\{k \in \N_0 \mid j+kd \in X\}\:.$$
\end{definition}

Clearly, we have
\begin{lemma}\label{d-corr} The map $s_d$
is a bijection between the set of all $\beta$-sets and
the set of all $d$-symbols.
\end{lemma}

\begin{definition}\label{asspart}
To an arbitrary $d$-symbol $S$,  we associate a partition
$p(S)$,  defined by
$$p(S)=p(s_d^{-1}(S)).$$
\end{definition}

\smallskip

\bigskip

\begin{definition}\label{hooks}
{\rm (1)} A {\rm hook} in  the $\beta$-set $X$ is a pair $(a,b)$ of
nonnegative integers with $a>b$
such that $a\in X$ and $b\notin X$.  The set of hooks in $X$ is denoted
$H(X)$.
If $a-b=\ell$ then $(a,b)$ is called an {\em $\ell$-hook} in $X$.

\smallskip

{\rm(2)} A {\rm hook} in the $d$-symbol   $S=(X_0,X_1,\ldots,X_{d-1})$
is a quadruple $(a,b,i,j)$ of nonnegative integers where
$a \ge b$, $i,j \in [d]$, $a \in X_i$, $b \notin X_j$,
and in addition if $a=b$ then $i>j$.
If $a=b$, we call the hook {\rm short} and otherwise {\rm long}.
The set of hooks in $S$ is denoted by~$H(S)$.
If $a-b=\ell$ and $(i-j)_{[d]} = e$, $e \in [d]$ then
$(a,b,i,j)$ is called an {\em $(\ell,e)$-hook} in~$S$.
\end{definition}

\smallskip

\noindent {\it Note.}
(i) To avoid confusion we want to point out that since we will be
dealing with generalized hook lengths below,
{\it an $\ell$-hook will {\rm not} be the same as a hook of length~$\ell$.}

\smallskip

(ii)  In \cite{Ma} only long hooks are considered in
symbols, since short hooks give only trivial contributions to the
unipotent degrees.
The inclusion of short hooks makes our later arguments much simpler,
and the short hooks will then be dealt with
separately whenever necessary.

\smallskip

\begin{lemma}\label{canbijs} 
{\rm (1)} There is a canonical bijection $\mathfrak{h}_X$ between the
  set $H(X)$ of hooks in a $\beta$-set $X$ and the set $H(p(X))$
  of hooks in the partition $p(X)$.  Thereby an $\ell$-hook in $H(X)$ is
  mapped to an $\ell$-hook in the partition $H(p(X))$.

\smallskip

{\rm (2)}  There is a canonical bijection $\mathfrak{h}_{X,d}$
between the  set $H(X)$ of hooks in a $\beta$-set $X$ and the
set $H(S)$ of hooks in the associated $d$-symbol $S=s_d(X)$.

\smallskip

{\rm (3)} There is a canonical bijection $\mathfrak{h}_S$ between
the set $H(S)$ of hooks in a $d$-symbol $S$
and the set $H(p(S))$ of hooks in the partition $p(S)$.
\end{lemma}

\noindent {\bf Proof.} Let $X=\{a_1,a_2,\ldots,a_t\}$.
(1) If $(a,b)\in  H(X)$, $a-b=\ell$ and $a=a_i$, then
the partition $p(X)$ has an $\ell$-hook in the $i$th row.
The map $\mathfrak{h}_X$ sending $(a,b)$ to this hook is a bijection;
this is well known, see \cite[Section 2.7]{JK}  or \cite[section 1]{OlL}
for more details.

\smallskip

(2) The map $\mathfrak{h}_{X,d}$ which maps $(a,b) \in H(X)$ to
$(a',b',i,j) \in H(S)$ where $a=a'd+i$, $b=b'd+j$ with $i,j \in [d]$,
is obviously a bijection between $H(X)$ and~$H(S)=H(s_d(X))$.
\smallskip

(3) By Definition~\ref{asspart}, $p(S)=p(s_d^{-1}(S))$. 
Put $\mathfrak{h}_S=\mathfrak{h}_X \circ\mathfrak{h}_{X,d}^{-1}$, where
$X=s_d^{-1}(S)$.  \qed

\bigskip
\begin{definition} \label{hookrem}
{\rm (1)} Let $X$ be a $\beta$-set and $z=(a,b) \in H(X)$. If
$X'=(X\setminus\{a\})\cup \{b\}$ we say
that $X'$ is obtained by {\rm removing the hook} $z$ from $X$.
\smallskip

{\rm (2)} Let $S=(X_0,\ldots,X_{d-1})$ be a $d$-symbol
and $z=(a,b,i,j) \in H(S)$.
If $i\neq j$, we set $X'_i=X_i\setminus \{a\}, X'_j=X_j \cup \{b\}$,
for $i=j$ we set $X'_i=(X_i\setminus \{a\})\cup \{b\}$,
and we set $X'_k=X_k$ for all $k \neq i,j$.
Then we say that $S'=(X'_0,\ldots,X'_{d-1})$ is
obtained by {\rm removing the hook} $z$ from $S$.
\end{definition}

\begin{rem} \label{cores}{\rm If we keep removing $\ell$-hooks from a $\beta$-set
  $X$ for a fixed $\ell$  we eventually reach a $\beta$-set with no
  $\ell$-hooks left.
  This is the $\ell$-core of~$X$, denoted $C_{\ell}(X)$. By an abacus
  argument this is well-defined \cite[ 2.7.16]{JK}.
  A similar statement is true for $(\ell,e)$-hooks in a
  symbol \cite[3.4]{Ma}. We return to this in Section~\ref{symbs}.}
\end{rem}

We call a $d$-symbol $S=(X_0,\ldots, X_{d-1})$  {\it balanced} if
$|X_0|=|X_1|= \cdots = |X_{d-1}|$, and if in addition there is
an $i\in [d]$ such that $0 \notin X_i$.

\smallskip

To a $d$-tuple $(\kappa_0,\kappa_1,\ldots,\kappa_{d-1})$ of partitions
we associate a balanced $d$-symbol as follows.
If $r$ is the maximal length (i.e., number of parts) of the partitions
$\kappa_i$, we may choose a $\beta$-set $Y_i$ of cardinality
$r$ for each partition~$\kappa_i$.  Then $0 \notin Y_j$ whenever
$\kappa_j$ has length $r$ because then $Y_j$ is just the set of first
column hook lengths of $\kappa_j$. Thus
$$t_d(\kappa_0,\kappa_1,\ldots,\kappa_{d-1}) =(Y_0,Y_1,\ldots, Y_{d-1})$$
is a well-defined balanced $d$-symbol for the $d$-tuple of partitions.
On the other hand, if $S=(X_0,\ldots, X_{d-1})$ is balanced,
then  $t_d(p(X_0),\ldots,p(X_{d-1}))=(X_0,\ldots, X_{d-1})$.

We have shown:
\begin{lemma}\label{symcorequo}
There is a bijection $t_d$ between the set of $d$-tuples of
partitions and the set of balanced $d$-symbols.
\end{lemma}

\begin{definition} \label{squotcore}
{\rm (1)} The  {\em balanced quotient} $Q(S)$ of an  arbitrary $d$-symbol
$S=(X_0,X_1,\ldots, X_{d-1})$ is
defined as the balanced $d$-symbol
$$Q(S)=t_d(p(X_0),p(X_1),\ldots,p(X_{d-1})).$$
We call $q(S)=p(Q(S))$ the {\em quotient partition of $S$.}
\smallskip

\noindent {\rm (2)} The  {\rm core} $C(S)$ of an  arbitrary $d$-symbol
$S=(X_0,X_1,\ldots, X_{d-1})$ is defined as the $d$-symbol
$$C(S)=([x_0],[x_1],\ldots,[x_{d-1}]),$$
where $x_i=|X_i|$ for $i \in [d]$. We call $c(S)=p(C(S))$ the {\rm core partition of $S$.}
\end{definition}

\begin{rem} {\rm A remark on notation. The core  $C(S)$ of a
    $d$-symbol $S$ is really its $(1,0)$-core (see Remark~\ref{cores}) and we
    will consider $Q(S)$ as the $(1,0)$-quotient of $S$. This will be
    generalized in Section~\ref{symbs}.}
\end{rem}

\begin{rem} {\rm We may recover a $d$-symbol $S=(X_0,X_1,\ldots,X_{d-1})$
  from its balanced quotient $Q(S)=(Y_0,Y_1,\ldots,Y_{d-1})$
  together with its core $C(S)=(Z_0,Z_1,\ldots,Z_{d-1})$. Indeed, for $i \in
  [d]$ $X_i$ must be the $\beta$-set of cardinality $|Z_i|$ for the
  partition $p(Y_i).$}
\end{rem}

\begin{definition}\label{Xdquot} 
Let $X$ be a $\beta$-set. Then the {\em $d$-quotient
  partition of $X$} is defined as $q_d(X)=q(s_d(X))$ and the
{\em $d$-core partition of $X$} is defined as $c_d(X)=c(s_d(X))$.

With notation as in Definition~\ref{asssymb}
we define the {\em $d$-quotient of $X$} as
$$Q_d(X)=s_d^{-1}(t_d(p(X^{(d)}_0),p(X^{(d)}_1),\ldots,p(X^{(d)}_{d-1})))\:.$$
\end{definition}

\begin{rem} \label{corequorem} {\rm If we put the elements of $X$ on the
$d$-abacus and if $X^{(d)}_j$ is as in Definition~\ref{asssymb} then the
results of  \cite[Section 2.7]{JK}  or \cite[section 1]{OlL} show
the following}

$\bullet$ ~ {\rm The  $d$-core partition $c_d(X)$ of $X$ is also the  
{\em $d$-core $\la_{(d)}$} of the partition $\la=p(X)$ and thus
it stays the same when we replace $X$ by $X^{+s}$.

$\bullet$ ~ We have $c_d(X)=p(C_d(X))$.

$\bullet$  ~ The  $d$-quotient partition $q_d(X)$ of $X$ must have
an} empty $d$-core. {\rm It  may}
change if we replace $X$ by $X^{+1}$. {\rm However $X$ and $X^{+d}$
have the same $d$-quotient partition.

$\bullet$  ~ We have $q_d(X)=p(Q_d(X))$ and $s_d(Q_d(X))=Q(s_d(X))$.

$\bullet$ If $t \in [d]$  is chosen such that $|X|+t$ is divisible by $d$ then
$$(p(X^{(d)}_t),p(X^{(d)}_{t+1}) \ldots,p(X^{(d)}_{t+d-1}))$$
is the}  $d$-quotient of the partition $p(X)$ \cite[2.7.29]{JK}.
{\rm  Here the subscripts are to be read modulo~$d$.}
\end{rem}

\medskip

\begin{lemma}\label{corequo} Let $X$ be a $\beta$-set, $d \in \N$. Then
$$|p(X)|=|q_d(X)|+|c_d(X)|.$$
\end{lemma}

\noindent {\bf Proof.} Using Remark~\ref{corequorem} we get
$$|p(X)|=d(|p(X^{(d)}_0)|+|p(X^{(d)}_1)|+\cdots+|p(X^{(d)}_{d-1})|)+|c_d(X)|.$$
If we apply this result to the $d$-quotient of $X$,
$$Q_d(X)=s_d^{-1}(t_d(p(X^{(d)}_0),p(X^{(d)}_1),\ldots,p(X^{(d)}_{d-1}))),$$
and the partition $q_d(X)=p(Q_d(X))$ (which has an empty $d$-core) we get
$|q_d(X)|=|p(Q_d(X))|=d(|p(X^{(d)}_0)|+|p(X^{(d)}_1)|+\ldots+|p(X^{(d)}_{d-1})|)$. 
The result follows. \qed

\medskip

The following diagram illustrates the connection between a $\beta$-set
$X$, its associated $d$-symbol $S=s_d(X)$, their cores and
quotients and the associated partitions via the partition map $p$:

\begin{equation*}
\xymatrix{
q_d(X) & Q_d(X) \ar[l]_p \ar[r]^{s_d}
  &Q(S) \ar[r]^-p& q(S)=q_d(X) \\
p(X)\ar[u]^{quot}\ar[d]_{core} &  {\it  X}  \ar[l]_p \ar[r]^{s_d}
\ar[u]^{quot}\ar[d]_{core} & {\it S} \ar[u]^{quot}\ar[d]_{core} \ar[r]^-p& p(S)=p(X) \ar[u]^{quot}\ar[d]_{core}\\
c_d(X) & C_d(X) \ar[l]_p \ar[r]^{s_d}
  &C(S)\ar[r]^-p& c(S)=c_d(X)
}
\end{equation*}

\medskip

Here is an example:

\begin{example} \label{ex1} {\rm
Consider the $\beta$-set $X=\{11,8,6,2,0\}$ for the partition
$p(X)=\la=(7,5,4,1)$ of 17. Let $d=3$.
%\iffalse{===========
The 3-abacus representation for $X$ is
\begin{center}
\begin{tabular}{c c c}
{\bf 0}& 1&{\bf 2} \\
3&4&5\\
{\bf 6}&7&{\bf 8} \\
9&10&{\bf 11} \\
\end{tabular}
\end{center}
%=========}\fi We compute
corresponding to $S=s_3(X)=(\{2,0\},\emptyset,\{3,2,0\})$
(giving the levels of the beads).
We get the balanced quotient
$$Q(S)=t_3((1),(0),(1,1))= (\{2,0\},\{1,0\},\{2,1\}).$$
Below on the left, we see its  3-abacus representation
and on the right, the 3-abacus representation of the
$\beta$-set corresponding to $Q(S)$
\begin{center}
$Q(S)$ \quad \begin{tabular}{c c c}
{\bf 0} &  {\bf 0} & 0 \\
1 & {\bf 1} & {\bf 1} \\
{\bf 2} & 2 & {\bf 2} \\
3&3&3\\
\end{tabular}
\qquad
$Q_3(X)$ \quad \begin{tabular}{c c c}
{\bf 0} &  {\bf 1} & 2 \\
3 & {\bf 4} & {\bf 5} \\
{\bf 6} & 7 & {\bf 8} \\
9&10&11\\
\end{tabular}
\end{center}
That is,
$$Q_3(X)=s_3^{-1}(Q(S))=\{8,6,5,4,1,0\},$$
giving the 3-quotient partition
$$q_3(X)=(3,2,2,2).$$
The cardinalities of the $\beta$-sets in $S$ give the core symbol
$$C(S)=((\{1,0\},\emptyset,\{2,1,0\}).$$
The corresponding $\beta$-set (and by construction also the 3-core of $X$)
is
$$C_3(X)=\{8,5,3,2,0\},$$
and thus the corresponding 3-core partition
is $$c_3(X)=(4,2,1,1).$$
Note that $|q_3(X)|+|c_3(X)|= 9+8=17=|p(X)|$, illustrating the
previous lemma.}
\end{example}

\begin{lemma}\label{card} Let $X$ be a $\beta$-set
%for the partition  $p(X)$
and $S=s_d(X)$ be the associated $d$-symbol.
%Let $Q_d(X)$ and $C_d(X)$ be
%$\beta$-sets with $s_d(Q_d(X))=Q(S)$ and  $s_d(C_d(X))=C(S)$.
We have:
$$
\begin{array}{ccccc}
 (1) & |p(X)|&=|H(p(X))|&=|H(X)|&=|H(S)| \\
 (2) &|q_d(X)|&= |H(q_d(X))|&=|H(Q_d(X))|&=|H(Q(S))| \\
 (3) &|c_d(X)|&=|H(c_d(X))|&=|H(C_d(X))|&=|H(C(S))| \\
\end{array}
$$
and in addition
$$
\begin{array}{lllll}
 (4)& |H(S)|=|H(Q(S))|+|H(C(S))|&&&
\end{array}
$$

\end{lemma}

\noindent {\bf Proof.} (1)-(3) are trivial using the bijections in
Lemma~\ref{canbijs} and then (4) follows from Lemma~\ref{corequo}. \qed

%%%%%%%%%%%%%%%%\noindent {\bf 3. Correspondences between sets of hooks.}
\section{Correspondences between sets of hooks} \label{hookscor}

In this section we fix the following notation.
Let $S=(X_0,X_1,\ldots,X_{d-1})$ be a $d$-symbol
with balanced quotient symbol
$Q=Q(S)=(Y_0,Y_1,\ldots,Y_{d-1})$
and core symbol
$C=C(S)=(Z_0,Z_1,\ldots,Z_{d-1})$.\\
We would like to describe some well-behaved correspondences (called
{\it universal bijections})  between the set
$ H(S)$ and the  union of the sets $H(Q)$ and $ H(C)$ (see Lemma~\ref{card}(4)).

We decompose the set of hooks $H(S)$ into {\it disjoint} subsets
according to the position of the hooks. For this, define
$$H_{ij}(S)=\{(a,b,i,j) \mid (a,b,i,j) \in H(S)\}\, ,
\: H_i(S) = H_{ii}(S)\:.$$
Then we have
$$ H(S)=(\bigcup_{i \in [d]} H_i(S)) \cup
(\bigcup_{i,j \in [d],~i \neq j}H_{ij}(S))\:.$$
Also, we set $H_{\{ij\}}(S)=
H_{ij}(S) \cup H_{ji}(S)$.

We want to split  $H_{ij}(S)$ further according to the
differences $a-b$,  i.e., for $\ell \ge 0$ we define
$$ H_{ij}^{\ell}(S)=\{(a,b,i,j) \in H_{ij}(S)~|~a-b=\ell\}\:, \:
H_i^{\ell}(S)=H_{ii}^{\ell}(S)\:.$$

We have the following easy observation:
\begin{lemma}\label{card2}
Let $i,j \in [d]$. Then we have
$$
|H_{ij}^{\ell}(S)|  =
\left\{
\begin{array}{cl}
|X_i|-|X_i \cap X_j^{+\ell}|  & \text{if $\ell >0$, or $\ell=0$ and  $i>j$}\\
0 & \text{if $\ell=0$ and $i \le j$}
\end{array} \right.
 \:.$$
\end{lemma}

Towards our key result,
we proceed  to describe a correspondence between
 the set  $H_{\{ij\}}(S)$  and the multiset union
of the sets  $H_{\{ij\}}(Q)$ and $H_{\{ij\}}(C)$,
for any given $i,j$.

Note that for each $i \in [d]$ and each $\ell \ge 0$
there is a bijection $H_i^{\ell}(S) \to H_i^{\ell}(Q)$
because $p(X_i)=p(Y_i)$, and $H_i^{\ell}(C)=\emptyset$.
\\
Thus it suffices to consider the situation where $i \ne j$.
We may assume that $i$ and $j$ are such that  $\Delta=|X_i|-|X_j| \ge 0$.

For $\Delta=0$ we clearly have $H_{\{ij\}}(S)=H_{\{ij\}}(Q)$ and $H_{\{ij\}}(C)=\emptyset$.
\\
So now we consider the case $\Delta>0$.
We may  assume that
$X_i=Y_i$, $Y_j=X_j^{+\Delta}$, $Z_i=[\Delta]$ and
$Z_j=\emptyset$.

Note that for $1 \le \ell \le \Delta$,
$H_{ij}^{\ell}(C)$ has cardinality $\Delta-\ell$.
Also $H_{ij}^0(C)$ has cardinality 0 if $i<j$
and cardinality $\Delta$ if $j<i$.
Clearly, $H_{ij}^{\ell}(C)=\emptyset$ if $\ell > \Delta$,
and $H_{ji}^{\ell}(C)=\emptyset$.

First, let $\ell >\Delta $.
The conditions $X_i=Y_i$, $X_j^{+\Delta}=Y_j$
imply that we have bijections
$$H_{ji}^{\ell-\Delta}(S) \to H_{ji}^{\ell}(Q),\:
(a,b,j,i) \mapsto (a+\Delta,b,j,i)\:$$
$$H_{ij}^{\ell}(S) \to H_{ij}^{\ell-\Delta}(Q), \:
(a,b,i,j) \mapsto (a,b+\Delta,i,j).$$
Next we consider  $H_{ji}^{\ell}(Q)$  and $H_{ij}^{\ell}(S)$
for $0 < \ell < \Delta$.
By Lemma~\ref{card2} we obtain
$$|H_{ji}^{\ell}(Q)|=
|X_j|+\Delta -|X_j^{+\Delta}\cap X_i^{+\ell}|
=|X_j|+(\Delta-\ell)-|X_j^{+(\Delta-\ell)}\cap X_i|$$
and then
$$|H_{ij}^{\Delta-\ell}(S)|=|X_i|-|X_i\cap X_j^{+(\Delta-\ell)}|
=|H_{ji}^{\ell}(Q)| +\ell \:.$$
Replacing $\ell$ by $\Delta-\ell$ we have for $0 < \ell < \Delta$:
$$|H_{ij}^{\ell}(S)|=|H_{ji}^{\Delta-\ell}(Q)|+|H_{ij}^{\ell}(C)|\:.$$
It remains to consider  the contributions of the sets
$H_{ij}^{\Delta}(S)$, $H_{ij}^{0}(S)$, $H_{ji}^{0}(S)$ towards
$H_{\{ij\}}(S)$, and of
$H_{ji}^{\Delta}(Q)$, $H_{ji}^{0}(Q)$, $H_{ij}^{0}(Q)$
towards $H_{\{ij\}}(Q)$, and take the contribution
from $H_{ij}^0(C)$ into account.
We use again Lemma~\ref{card2} and keep in mind that
$X_j^{+\Delta}=Y_j$.
Thus we have

$|H_{ij}^{\Delta}(S)|=|X_i|-|X_i \cap Y_j|$,
$|H_{ji}^{\Delta}(Q)|=|Y_j|-|Y_j \cap Y_i^{+\Delta}|
=|X_j|-|X_j\cap Y_i|$.
\\[5pt]
Furthermore, for $i>j$ we have

$|H_{ij}^{0}(S)|= |X_i|-|X_i \cap X_j|$,
$|H_{ji}^{0}(S)|=  0$,

\smallskip

$|H_{ji}^{0}(Q)|=0$,
$|H_{ij}^{0}(Q)|=|Y_i|-|Y_i \cap Y_j|$,
$|H_{ij}^0(C)|=\Delta$,
\\[5pt]
and for  $i<j$ we have

$|H_{ij}^{0}(S)|= 0$,
$|H_{ji}^{0}(S)|= |X_j|-|X_j \cap X_i|$,

\smallskip

$|H_{ji}^{0}(Q)|=|Y_j|-|Y_j \cap Y_i|$,
$|H_{ij}^{0}(Q)|=0$,
$|H_{ij}^0(C)|=0$.
\\[5pt]
Let $k=\max(i,j)$;
using that $X_i=Y_i$ and $|Y_j|=|Y_i|=|X_j|+\Delta$,
we see that in both cases the terms add up:
 $$\begin{array}{rcl}
 |H_{ij}^{\Delta}(S)|+|H_{\{ij\}}^{0}(S)|
 &=&
|X_i|-|X_i\cap Y_j|+|X_k|-|X_i\cap X_j|\\[8pt]
&=&
|X_j| - |Y_i\cap X_j| + |X_k|+\Delta - |Y_i \cap Y_j| \\[8pt]
&=&
|H_{ji}^{\Delta}(Q)|+|H_{\{ij\}}^{0}(Q)|+|H_{ij}^0(C)|.
\end{array}$$
\medskip

We have shown the following key result:
\begin{prop} \label{Prop:key} Let  $S=(X_0,\ldots,X_{d-1})$ be
a  $d$-symbol with balanced quotient symbol $Q(S)=Q$ and core symbol $C(S)=C$.
Let $i \ne j \in [d]$, and set $\Delta =|X_i|-|X_j|$.
\\
When $\Delta > 0 $,
we have the following equalities:

$\bullet$ For all $\ell > \Delta$:  $|H_{ij}^{\ell}(S)|=|H_{ij}^{\ell-\Delta}(Q)|$.

\smallskip

$\bullet$  For all $\ell > \Delta$: $|H_{ji}^{\ell-\Delta}(S)|=|H_{ji}^{\ell}(Q)|$.

\smallskip

$\bullet$  For all $0 <\ell < \Delta $:
$|H_{ij}^{\ell}(S)|=|H_{ji}^{\Delta-\ell}(Q)|+|H_{ij}^{\ell}(C)|$.

\smallskip

$\bullet$ For $\ell = \Delta $: $|H_{ij}^{\Delta}(S)|=
\left\{\begin{array}{cl}
|H_{ij}^{0}(Q)| = |H_{\{ij\}}^{0}(Q)|  & \text{if } i> j\\
|H_{ji}^{0}(Q)| = |H_{\{ij\}}^{0}(Q)| & \text{if } i< j\\
\end{array}\right. $.

\smallskip

$\bullet$ For $\ell = 0 $:
$|H_{ji}^{\Delta}(Q)|+|H_{ij}^0(C)|=
\left\{\begin{array}{cl}
|H_{ij}^{0}(S)| = |H_{\{ij\}}^{0}(S)|  & \text{if } i> j\\
|H_{ji}^{0}(S)| = |H_{\{ij\}}^{0}(S)| & \text{if } i< j\\
\end{array}\right. $.

\smallskip

$\bullet$  $|H_{ij}^{\Delta}(S)|+|H_{\{ij\}}^{0}(S)|=
|H_{ji}^{\Delta}(Q)|+| H_{\{ij\}}^{0}(Q)|+|H_{ij}^0(C)|$.
\\[5pt]
When $\Delta=0$, we have

\smallskip

$\bullet$ $|H_{ij}^{\ell}(S)|=|H_{ij}^{\ell}(Q)|$, $H_{ij}^{\ell}(C)=\emptyset$,
for all $\ell \ge 0 $.
\end{prop}

Using all the correspondences behind the
equalities established so far
%of Proposition~\ref{Prop:key} for all relevant
%$i,j$ together with the correspondences between $H_i^{\ell}(S)$  and
%$H_i^{\ell}(Q)$
we find bijective correspondences (between multisets)
$$ H_{\{ij\}}(S) \to H_{\{ij\}}(Q) \cup H_{\{ij\}}(C)$$
for all $i,j\in [d]$, and we may glue these together to set up
a {\it universal bijection}
$$\omega_S: H(S) \to H(Q) \cup H(C)\:.$$

\begin{rem}\label{ellpos}
{\rm
Let us specify the properties of a universal bijection $\omega_S$
very explicitly.
Let $i \neq j\in [d]$ be chosen such that $\Delta=|X_i|-|X_j| \ge 0$.
\\
(i) For positive $\ell$ we have:

 {\rm (1)} For all $\ell > \Delta \ge 0$:  $\omega_S(H_{ij}^{\ell}(S))=H_{ij}^{\ell-\Delta}(Q)$.

 {\rm (2)} For all $\ell > \Delta \ge 0$: $\omega_S(H_{ji}^{\ell-\Delta}(S))=H_{ji}^{\ell}(Q)$.

 (3) For $\ell = \Delta >0$: $\omega_S(H_{ij}^{\Delta}(S))=
\left\{\begin{array}{cl}
H_{ij}^{0}(Q) & \text{if } i> j\\
H_{ji}^{0}(Q) & \text{if } i< j\\
\end{array}\right. $

(4) For all $0 < \ell < \Delta$:
$\omega_S(H_{ij}^{\ell}(S))=H_{ji}^{\Delta-\ell}(Q)\cup
H_{ij}^{\ell}(C)$.

(5) For all $\ell$: $\omega_S(H_i^{\ell}(S))=H_i^{\ell}(Q)$.
\\[1ex]
(ii) \label{ello}  For  $\ell=0$, we have the following.

For $\Delta>0$:

(1)
$
H_{ji}^{\Delta}(Q) \cup H_{ij}^0 (C) =
\left\{ \begin{array}{cl}
\omega_S(H_{ij}^{0}(S))=
 \omega_S(H_{\{ij\}}^{0}(S)) & \text{if } i> j\\
 \omega_S(H_{ji}^{0}(S))=
 \omega_S(H_{\{ij\}}^{0}(S)) & \text{if } i< j\
\end{array}\right. $

%{\rm (1)} If\ $i>j$, then $\omega_S(H_{ij}^0(S)) = H_{ji}^{\Delta}(Q) \cupH_{ij}^0(C)$.

%{\rm (2)} If $i<j$ then  $\omega_S(H_{ji}^0(S))= H_{ij}^{\Delta}(Q)$.

{\rm (2)} $\omega_S(H_i^{0}(S))=H_i^{0}(Q)$.

For $\Delta=0$:

{\rm (3)} $\omega_S(H_{ij}^{0}(S))=H_{ij}^{0}(Q)$. }
\end{rem}

Of course such a universal bijection is by no means unique and
apparently cannot be made ``canonical''. The important fact for this
bijection is the relation between the   $\ell$'s in corresponding
hooks.

\bigskip

%\noindent {\bf 4. Hook length functions}
\section{Generalized hook length functions}\label{hookfunc}

\medskip
We now want to associate lengths to the hooks in $H(S)$,   where
$S$ is a $d$-symbol.
We define our (generalized) hook length function on the set
$$H = \{(a,b,i,j)\in \N_0^2\times [d]^2 \mid a\ge b ~~{\rm and}~~ \, i>j \text{ if } a=b\}\:.$$
In general the lengths may be arbitrary real numbers, i.e., we have
a (generalized) hook length function $h: H \to \R$.
However, we only want to consider
functions $h$ such that

{\it the value $h(a,b,i,j)$  depends only on $\ell=a-b$, $i$ and~$j$.}
\\
This guarantees that
the multiset $\mathcal{H}(S)$ of all
$h(z)$, $z \in H(S)$,  coincides with
$\mathcal{H}(S^{+s})$ for all $s \in \N_0$.
Indeed, $(a,b,i,j) \in
H(S)$ if and only if $(a+s,b+s,i,j) \in
H(S^{+s})$, and the $h$-value for these hooks will be the same.
\\
We set
$$H_{ij} = \{(a,b,i,j) \mid (a,b,i,j) \in H\} \, , \:
H_{ij}^{\ell} = \{(a,b,i,j) \in H_{ij} \mid a-b=\ell\}$$
and then the hook length functions $h$ that we will consider will
be constant on $H_{ij}^{\ell}$.
\smallskip

We now describe some hook length functions of interest for $d$-symbols.
For a $(d+1)$-tuple  $\delta=(c_0,c_1,\ldots,c_{d-1};k)$
of real numbers, with  $k\ge 0$,
we define the {\em $\delta$-length}
of  $(a,b,i,j)\in H$  as
$$ h^{\delta}(a,b,i,j) = k(a-b)+c_i-c_j\:.$$
We call $\delta$ a {\it d-hook data tuple}.
For any $d$-symbol $S$,
we let $\mathcal{H}^{\delta}(S)$ be the multiset of all
$h^{\delta}(a,b,i,j)$, $(a,b,i,j) \in H(S)$,
and  $\mathcal{H}^{\delta}_{ij}(S)$ be the multiset of all
$h^{\delta}(a,b,i,j)$, $(a,b,i,j) \in H_{ij}(S)$.

\begin{rem}\label{choice}
{\rm
Some special choices of $d$-hook data tuples will be particularly
important in the next sections.

\smallskip

$\bullet$ If we choose $\delta^o=(0,0,\ldots, 0;1)$ then the
$\delta$-length of long hooks in
$S$ coincides with the length defined in \cite[p.~782]{Ma},  and the
short hooks have $\delta^o$-length~0. We call  $\delta^o=(0,0,\ldots,
0;1)$ the {\it minimal d-hook data tuple}.

\smallskip

$\bullet$ If we choose  $\delta^*=(0,1,\ldots, d-1;d)$ then in the
notation of Lemma~\ref{canbijs} the usual hook length $a-b$ of $(a,b)$
in $H(X)$ equals the $\delta^*$-length of the corresponding
hook $\mathfrak{h}_{X,d}(a,b)$ in $H(S)$.
We call  $\delta^*=(0,1,\ldots,d-1;d)$ the {\it partition d-hook
  data tuple}.
}
\end{rem}

\smallskip

As before,
we let
$S=(X_0,X_1,\ldots,X_{d-1})$
be a $d$-symbol with balanced quotient
$$Q=Q(S)=(Y_0,Y_1,\ldots,Y_{d-1})$$ and core
$$C=C(S)=([x_0],[x_1],\ldots,[x_{d-1}]),$$
where $x_i=|X_i|$.
We have set up above a universal bijection
$$\omega_S: H(S) \rightarrow H(Q) \cup H(C),$$
with properties specified in Remark~\ref{ellpos}.

Let $\delta=(c_0,c_1,\ldots,c_{d-1};k)$ be an arbitrary $d$-hook data tuple.
We define
$$\delta_S:=(c_0+x_0k,c_1+x_1k,\ldots ,c_{d-1}+x_{d-1}k;k).$$
This is a new $d$-hook data tuple which depends on the core of~$S$.

We want to modify $h=h^{\delta}$ to a new length function
$\overline{h}=\overline{h}^{\delta_S}$
such that for $z\in H(S)$ we have

\bdm
h^{\delta}(z)=\left\{ \begin{array}{lcl}
\overline{h}^{\delta_S}(\omega_S(z)) &\textrm{if}& \omega_S(z) \in H(Q) \\
h^{\delta}(\omega_S(z)) &\textrm{if }& \omega_S(z) \in H(C).
\end{array} \right.
\edm

If this is done we immediately have
$$\mathcal{H}^{\delta}(S)=\overline{\mathcal{H}}^{\delta_S}(Q)
 \cup \mathcal{H}^{\delta}(C), $$
where $\overline{\mathcal{H}}^{\delta_S}(Q)$ is the multiset of all
$\overline{h}^{\delta_S}(z)$, $z \in H(Q)$.

We proceed to define $\overline{h}^{\delta_S}$. Apart from a sign,
$\overline{h}^{\delta_S}$ is just $h^{\delta_S}$,
 with $\delta_S$ as above. Let us describe the sign modification of  $h^{\delta_S}$ on
$H_{\{ij\}}=H_{ij}\cup H_{ji}$. \\
We set $H_{ij}^{>m} = \bigcup_{\ell >m} H_{ij}^{\ell}$,
and use similar notation for the condition $\ge m$ and
for subsets of $H_{ji}$.
\\
We assume that $i,j\in [d]$ are such that $\Delta=x_i-x_j\ge 0$.
Then for $z \in  H_{\{ij\}}^{\ell}$ we define
\bdm
\overline{h}^{\delta_S}(z) =\left\{ \begin{array}{rl}
h^{\delta_S}(z) & \textrm{if } z \in H_{ij}^{\geq 0} \cup H_{ji}^{>\Delta},
\text{ or } z\in H_{ji}^{\Delta} \text{ if } i<j\\
-h^{\delta_S}(z)& \textrm{otherwise}
\end{array} \right.
\edm

\begin{theorem}\label{key}
Let $\omega_S$ be a universal bijection.
Let $h=h^{\delta}$ and $\overline{h}=\overline{h}^{\delta_S}$
be as defined above.  Then for $z \in H(S)$ we have
\bdm
h^{\delta}(z)=
\left\{ \begin{array}{lcl}
\overline{h}^{\delta_S}(\omega_S(z)) & \textrm{if}& \omega_S(z) \in H(Q)\\
h^{\delta}(\omega_S(z))&  \textrm{if} &  \omega_S(z) \in H(C)\\
\end{array} \right.
\edm
\end{theorem}
\smallskip
\noindent {\bf Proof.} Let $i,j \in [d]$.
As above, we assume that $i,j$ are chosen such that $\Delta:=x_i-x_j \ge 0 $.
For $i=j$, the claim clearly holds, so we assume now that $i \neq j$.
We  refer to Remark~\ref{ellpos} %and~\ref{ell0}
for the  properties of $\omega_S$ used below.

First assume $\Delta=0$.
Then $\omega_S(H^{\ell}_{ij}(S))=H^{\ell}_{ij}(Q)$
for all $\ell \ge 0$,
and by definition we have for $z \in H^{\ell}_{ij}(S)$
(and analogously for $z\in H_{ji}^{\ell}(S)$):
$$h^{\delta}(z)=\ell k+ c_i-c_j=h^{\delta_S}(\omega_S(z))=
\overline{h}^{\delta_S}(\omega_S(z))\:.$$

Now assume $\Delta>0$. First consider the case $\ell > \Delta$.
Take $z \in H_{ij}^{\ell}(S)$. Then  we know that
$\omega_S(z) \in H_{ij}^{\ell-\Delta}(Q)$ and by definition
$$
\begin{array}{rcl}
\overline{h}^{\delta_S}(\omega_S(z))&=&
h^{\delta_S}(\omega_S(z))= (\ell-\Delta)k+(c_i+x_ik)-(c_j+x_jk)\\[5pt]
&=& \ell k +c_i-c_j=h^{\delta}(z)\:.
\end{array}$$
Also, if $z \in H_{ji}^{\ell-\Delta}(S)$ then
$\omega_S(z) \in H_{ji}^{\ell}(Q)$ and again
$$
\begin{array}{rcl}
\overline{h}^{\delta_S}(\omega_S(z))&=&
h^{\delta_S}(\omega_S(z))= \ell k+(c_j+x_jk)-(c_i+x_ik)\\[5pt]
&=& (\ell-\Delta) k+c_j-c_i=h^{\delta}(z)\:.
\end{array}$$

Suppose next that $0< \ell < \Delta$.
If $z \in H_{ij}^{\ell}(S)$ then either
$\omega_S(z) \in H_{ji}^{\Delta-\ell}(Q)$ or
$\omega_S(z) \in H_{ij}^{\ell}(C)$.
In the latter case, clearly $h^{\delta}(\omega_S(z)) = h^{\delta}(z)$.
In the former case we compute
$$\begin{array}{rcl}
\overline{h}^{\delta_S}(\omega_S(z))&=&
-h^{\delta_S}(\omega_S(z))= (\ell-\Delta) k-(c_j+x_jk)+(c_i+x_ik)\\[5pt]
&=& \ell k -c_j+c_i=h^{\delta}(z)\:.
\end{array}
$$

If $\ell=0$,  $z \in H_{\{ij\}}^{0}(S)$,
we have to distinguish the cases $i>j$ and $i<j$.
When $i>j$,  $z\in H_{ij}^{0}(S)$, and $z$ can be mapped to either
$\omega_S(z) \in H_{ji}^{\Delta}(Q)$ or to
$\omega_S(z) \in H_{ij}^{0}(C)$.
The latter case is clear, and in the former case we compute
$$\overline{h}^{\delta_S}(\omega_S(z))= -h^{\delta_S}(\omega_S(z))
=-\Delta k - (c_j+x_jk)+(c_i+x_ik)=c_i-c_j=h^{\delta}(z)\:.
$$
When $i<j$,   $z\in H_{ji}^{0}(S)$, $z$ is mapped to
$\omega_S(z) \in H_{ji}^{\Delta}(Q)$,
and we compute
$$\overline{h}^{\delta_S}(\omega_S(z))= h^{\delta_S}(\omega_S(z))
=\Delta k + (c_j+x_jk)-(c_i+x_ik)=c_j-c_i=h^{\delta}(z)\:.
$$

Finally, we are in the case $z\in H_{ij}^{\Delta}(S)$.
Then $\omega_S(z)\in H_{\{ij\}}^0(Q)$, and again
we have to distinguish the cases $i>j$ and $i<j$.
When $i>j$, $\omega_S(z)\in H_{ij}^0(Q)$ and we have
$$\overline{h}^{\delta_S}(\omega_S(z))= h^{\delta_S}(\omega_S(z))
=c_i+x_ik-(c_j+x_jk)=\Delta k + c_i-c_j=h^{\delta}(z)\:.
$$
When $i<j$, $\omega_S(z)\in H_{ji}^0(Q)$ and we have
$$\overline{h}^{\delta_S}(\omega_S(z))= - h^{\delta_S}(\omega_S(z))
=- (c_j+x_jk)+c_i+x_ik =\Delta k + c_i-c_j=h^{\delta}(z)\:.
$$

Now we have dealt with all the elements in $H_{\{ij\}}(S)$ and
the assertion is proved.
 \qed
\medskip

As indicated before, the above theorem has the following consequence, which
we will use repeatedly in the following.

\begin{theorem}\label{metasatz} 
Let $S=(X_0,X_1,\ldots,X_{d-1})$ be a $d$-symbol with balanced quotient
$Q=Q(S)$ and core
$C=C(S)=([x_0],[x_1],\ldots,[x_{d-1}])$,
where $x_i=|X_i|$. Let $\delta=(c_0,c_1,\ldots,c_{d-1};k)$ be a $d$-hook data tuple
and $\delta_S=(c_0+x_0k,c_1+x_1k,\ldots ,c_{d-1}+x_{d-1}k;k)$. Then  we have the multiset equality
$$\mathcal{H}^{\delta}(S)=\overline{\mathcal{H}}^{\delta_S}(Q)
 \cup \mathcal{H}^{\delta}(C), $$
where $\overline{\mathcal{H}}^{\delta_S}(Q)$ is the multiset of all
$\overline{h}^{\delta_S}(z)$, $z \in H(Q)$.\\
In particular we have the multiset inclusion
 $$\mathcal{H}^{\delta}(C) \subseteq  \mathcal{H}^{\delta}(S)\:.$$
\end{theorem}

For a later application we need the following

\begin{rem}(Reversal of short hooks) \label{hookrev}
{\rm
We consider another possible definition of a short hook in a symbol,
which we call a ``reversed'' short hook. Suppose that
$S=(X_0,X_1,\ldots,X_{d-1})$. By Definition~\ref{hooks}
a short hook in $S$ is given by $(a,a,i,j)$
where $a \in X_i$, $a \notin X_j $ and $i >j$. Clearly short hooks from
$i$ to $j$ are determined by the elements
$a \in X_i \setminus (X_i\cap X_j)$.
The {\em reversed short hooks} in $S$ are  given by $(a,a,i,j)$
where $a \in X_i, a \notin X_j $ and  $i  < j$.
If $i>j$ and
$\delta=(c_1,\ldots,c_{d-1};k)$ is a $d$-hook data tuple then the
$\delta$-length of a short hook from $i$ to $j$ is $c_i-c_j$ whereas
the  $\delta$-length of a reversed short hook from $j$ to $i$ is
$c_j-c_i$. We note that these lengths are equal up to a sign.

Whenever $i >j$ then the number of short hooks from $i$ to $j$ is
 \hbox{$|X_i \setminus (X_i \cap X_j)|$}, whereas the number of reversed short
 hooks from $j$ to $i$ is $|X_j \setminus (X_i \cap X_j)|$. If the
 symbol $S$ is {\em balanced}, these numbers are equal, since then
 $|X_i|=|X_j|$.}
\end{rem}

%%%%%%%%%%%%%%%%%%%%%%%%%%%%%%%%%%%%%%%%%%%%%%%%%%%%%%%%%%%%%%%
\section{Partition data tuples and hooks in Partitions}
\label{appart}

\medskip

For any partition $\la$ we denote by $\mathcal{H}(\la)$ the
multiset of (usual) hook lengths in $\la$.
Let $X$ be any $\beta$-set and $d \in \N$. We let $\mathcal{H}(X)$ be
the multiset of hook lengths in $X$ (in the sense of
Definition~\ref{hooks}(1)).
Note that all these hook lengths are
positive integers. The bijection $\mathfrak{h}_X$ from Lemma~\ref{canbijs} (1)
preserves hook lengths. Thus we have

\begin{lemma}\label{betaparthook}  Let $X$ be a $\beta$-set for the
  partition $p(X)$. Then
$$\mathcal{H}(X)=\mathcal{H}(p(X)).$$
\end{lemma}

Consider the associated $d$-symbol
$S=s_d(X)$ to $X$. Lemma~\ref{canbijs} shows that there is a
bijection $\mathfrak{h}=\mathfrak{h}_{X,d}$ between
$H(X)$ and $H(S)$. If we choose the partition $d$-hook data tuple
$\delta^*=(0,1,\ldots,d-1;d)$ for $S$ then the description of $\mathfrak{h}$
shows that for any hook $z \in H(X)$ we have
$h(z)=h^{\delta^*}(\mathfrak{h}(z))$.

Thus we obtain
\begin{lemma}\label{betasymb}  Let $X$ be a $\beta$-set with
  associated $d$-symbol $S=s_d(X)$. Then
$$\mathcal{H}(X)=\mathcal{H}^{\delta^*}(S),$$
where $\delta^*$ is the partition $d$-hook data tuple.
\end{lemma}

If we then apply Theorem~\ref{metasatz} to a $d$-symbol $S$ and the
partition $d$-hook data tuple $\delta^*$ we get the following result:

\begin{theorem}\label{betameta}  Let the $d$-symbol $S=(X_0,X_1,\ldots,X_{d-1})$
have the balanced quotient $Q=Q(S)$
and the core $C=C(S)=([x_0],[x_1],\ldots,[x_{d-1}])$,
where $|X_i|=x_i$ for $0 \le i \le d-1$.  Let the $d$-hook data tuple
$\delta^*_S$ be defined by
$\delta^*_S=(x_0d, 1+x_1d,\ldots,(d-1)+x_{d-1}d;d)$. Then
$$\mathcal{H}^{\delta^*}(S)=\mathcal{H}^{\delta^*}(C) \cup
{\rm abs}(\mathcal{H}^{\delta^{*}_S}(Q))$$
where \; ${\rm abs}(\mathcal{H}^{\delta^{*}_S}(Q))=\{|h|~|~h \in
  \mathcal{H}^{\delta^{*}_S}(Q)\}. $
\end{theorem}

\noindent {\bf Proof.} In the notation of
Theorem~\ref{metasatz}, an element
$\overline{h}^{\delta^{*}_S}(z) \in \overline{\mathcal{H}}^{\delta^{*}_S}(Q)$,
$z \in H(Q)$,
has the same absolute value as
$h^{\delta^{*}_S}(z)$. Since the elements of
$\mathcal{H}^{\delta^*}(S)$ are positive integers,
the result follows.  \qed

\medskip

We may translate Theorem~\ref{betameta} into a statement
about hooks in partitions; 
we use the notation introduced in 
Definitions~\ref{asssymb} and~\ref{Xdquot}. 

\begin{theorem}\label{partmeta} 
Let $d \in \N$ and let $\la$ be a partition.
% with $d$-core $\la_{(d)}$.  
Let $X$ be a $\beta$-set for $\la$; we set $x_i=|X_i^{(d)}|$, for $i\in [d]$,   
and put $\delta=(x_0d, 1+x_1d,\ldots,(d-1)+x_{d-1}d;d)$. 
Let  $Q=Q(s_d(X))$, a balanced symbol. 
Then
$$\mathcal{H}(\la)=\mathcal{H}(\la_{(d)}) \cup
{\rm abs}(\mathcal{H}^{\delta}(Q))$$
where ${\rm abs}(\mathcal{H}^{\delta}(Q))=\{|h|~|~h \in
  \mathcal{H}^{\delta}(Q)\}. $
\end{theorem}

Let us add a remark to the theorem.
\begin{rem}
{\rm Theorem~\ref{partmeta} states in particular that the multiset of
hook lengths of the $d$-core
of a partition is contained in that of the partition. 
This suggests the following question. Suppose that $\mu$ is obtained 
from $\la$ by removing a number of $d$-hooks;  
is it then true that $\mathcal{H}(\mu) \subseteq \mathcal{H}(\la)?$
The answer is definitely {\em no}, 
and there are numerous examples for any~$d$. 
For instance, $\mu=(2d,1)$ is obtained  from
$\la=(2d,1^{d+1})$ by removing a $d$-hook,  
but $2d+1 \in \mathcal{H}(\mu)$ and 
$2d+1 \notin \mathcal{H}(\la)$.  }
\end{rem}

\smallskip

For any $d$-symbol $S$ there is a canonical bijection
$\mathfrak{h}_S$ between $H(S)$ and  $H(p(S))$ (Lemma~\ref{canbijs}).
This may be applied to the balanced symbol $Q$ 
of Theorem~\ref{partmeta}. 
Then $p(Q)$ is the $d$-quotient partition $q_d(X)$ of
the $\beta$-set $X$ of~$\la$. 
The lengths of corresponding
hooks differ up to a sign only by a multiple of~$d$. 
Thus $ \mathcal{H}^{\delta}(Q)$ may be seen as a multiset of
{\em modified hook lengths} of the $d$-quotient partition $q_d(X)$.

To be more specific we need the $d$-residues of the nodes in a Young
diagram. 
The node $(k,l)$ in row $k$ and column $l$ and the corresponding hook have
{\em ($d$-)residue} $e=(l-k)_{[d]} \in [d]$. 
%(\cite[2.7.34]{JK}). 
For a partition $\la=(\ell_1,\ldots, \ell_r)$, 
the rightmost residue in row $k$ is called the 
{\em hand ($d$-)residue} and the bottom residue in
column $l$ is called the {\em foot ($d$-)residue} 
of the $(k,l)$-hook in $\la$.
If $X=\{a_1,a_2,\ldots, a_s\}$ is a $\beta$-set
for $\la$, then $\ell_k=a_k-(s-k)$. 
Thus if $|X|=s$ is divisible by $d$, then
${a_k}_{[d]}=(\ell_k-k)_{[d]}$ 
is the end residue in row $k$ of $\la$. 
Note that if the $(k,l)$-hook has length~$m$, then 
its foot residue is congruent to $a_k-m+1$.  
Thus we have the following:

\begin{lemma}\label{handfoot} 
Let $X$ be a $\beta$-set for $\la$ such that $d \mid |X|$. 
Let $(a,b) \in H(X)$. Then the hand and foot ($d$-)residue of the
corresponding hook $\mathfrak{h}_X(a,b) \in H(\la)$ are $a_{[d]}$ and 
$b_{[d]}+1$,  respectively. 
\end{lemma}

\smallskip

We may then reformulate Theorem~\ref{partmeta}, involving 
a $d$-quotient partition instead of a balanced symbol, 
and a suitably modified hook length.

\begin{theorem} \label{partmeta2} 
Let $d\in \N$, $\la$ a partition, $X$ a $\beta$-set for $\la$,  
$x_i=|X_i^{(d)}|$, $i \in [d]$. 
%as  in Theorem~\ref{partmeta}. 
Let $\la^{(d)}_X=q_d(X)$ be the
  $d$-quotient partition of $X$. 
For $z \in H(\la^{(d)}_X)$, we define a modified hook length as 
$\overline{h}(z)=h(z)+(x_i-x_j)d$, if $z$ has hand and foot 
($d$-)residue $i$ and $j+1$ respectively. 
We denote by $\overline{\mathcal{H}}(\la^{(d)}_X)$
the multiset of all $\overline{h}(z)$, $z \in H(\la^{(d)}_X)$. 
Then
$$\mathcal{H}(\la)=\mathcal{H}(\la_{(d)}) \cup
{\rm abs}(\overline{\mathcal{H}}(\la^{(d)}_X))$$
where ${\rm abs}(\overline{\mathcal{H}}(\la^{(d)}_X))=\{|h|~|~h \in
  \overline{\mathcal{H}}(\la^{(d)}_X))\}$.
\end{theorem}

\noindent {\bf Proof.} 
It suffices to show that if $Q=Q(s_d(X))$ and $\delta$
are as in Theorem~\ref{partmeta}, then
$\mathcal{H}^{\delta}(Q)=\overline{\mathcal{H}}(\la^{(d)}_X)$. 
We have that $p(Q)=\la^{(d)}_X$ and that there is a canonical bijection
$\mathfrak{h}$ between $H(Q)$  and $H(\la^{(d)}_X)$ (Lemma~\ref{canbijs}).
In fact,  $\mathfrak{h}=\mathfrak{h}_{Q_d(X)}\circ
\mathfrak{h}_{Q_d(X),d}^{-1}$. 
Note that $Q=s_d(Q_d(X))$, and since $Q$ is balanced,  
we have $d \mid |Q_d(X)|$. 
Thus if $z \in H(\la^{(d)}_X)$ has hand and foot
($d$-)residues $i$ and $j+1$ respectively, then by Lemma~\ref{handfoot}
 $\mathfrak{h}_{Q_d(X)}^{-1}$ maps $z$ into a pair
$(a,b) \in H(Q_d(X))$ where  $a_{[d]}=i$ and $b_{[d]}= j$. 
Writing $a=a'd+i, b=b'd+j$ we have 
$\mathfrak{h}(z)=(a',b',i,j)$. 
Now $h(z)=a-b=(a'-b')d+i-j$,  
and by the definition of $\delta$ we have
$h^{\delta}(a',b',i,j)=\overline{h}(z)$.  
\qed

\bigskip

We illustrate Theorem~\ref{partmeta2} by an example. 
\begin{example} \label{ex2} (This continues Example~\ref{ex1}.) 
{\rm Consider the $\beta$-set $X=\{11,8,6,2,0\}$ for the partition
$p(X)=\la=(7,5,4,1)$ of 17. Let $d=3$. We computed
the 3-quotient partition $$\la^{(3)}_X=q_3(X)=(3,2,2,2)$$
and the 3-core partition
$$\la_{(3)}=c_3(X)=(4,2,1,1)\:.$$
Note that $|\la|=17=9+8=|\la^{(3)}_X|+|\la_{(3)}|$. 
Moreover, the numbers $x_i$ of elements in $X$ congruent to $i$ modulo
3 are (2,0,3).

Consider the hook diagrams of $\la$ and $\la_{(3)}$ where we have
marked by boldface eight hook lengths in $\la$ which also occur in $\la_{(3)}$.
\begin{center}
%{\scriptsize
$\la$ \quad \begin{tabular}{c c c c c c c}
$10$&$8$&${\bf 7}$&$6$&${\bf 4}$&${\bf 2}$&${\bf 1}$ \\
$7$&5&${\bf 4}$&$3$&${\bf 1}$ \\
$5$&$3$&${\bf 2}$&$1$& \\
${\bf 1}$ &&&& \\
\end{tabular}
%}\end{center}
\quad \quad 
%\begin{center}
%{\scriptsize
$\la_{(3)}$ \quad \begin{tabular}{c c c c}
$7$&$4$&$2$&$1$ \\
$4$&1& & \\
$2$&&& \\
$1$&&& \\
\end{tabular}
%}
\end{center}
The remaining nine hook lengths in $\la$ which are not in $\la_{(3)}$
are
$$R=\{1,3,3,5,5,6,7,8,10 \}.$$
We obtain these by adjusting the hook lengths of $\la^{(3)}_X$ by multiples
of 3 and changing signs of negative entries.
~Consider first the 3-residue diagram of~$\la^{(3)}_X$: 
\begin{center}
\begin{tabular}{c c c}
0&1&2 \\
2&0& \\
1&2& \\
0&1& \\
\end{tabular}
\end{center}

We add $3x_i$ to hook lengths in rows with end residue $i$ and
subtract $3x_j$ from hook lengths in columns with end residue $j+1$. 
These multiples of 3 are listed in boldface in the rows and columns of
the hook diagram of  $\la^{(3)}_X$ to the left, 
and we show the result (before sign change) to the right: 
\begin{center}
\begin{tabular}{l| c c c}
&${\bf 9}$&${\bf 6}$&${\bf 0}$ \\
\hline
${\bf 9}$&6&5&1 \\
${\bf 6}$&4&3& \\
${\bf 9}$&3&2& \\
${\bf 0}$&2&1& \\
\end{tabular}
\qquad \quad 
\begin{tabular}{r r r}
&&\\
$6$&$8$&$~10$ \\
$1$&3&  \\
$3$&$5$& \\
$-7$&$-5$& \\[15pt]
\end{tabular}
%}
\end{center} 
Changing the sign of -7 and -5 we get exactly the hook lengths
in the list ~$R$. }
\end{example}

\medskip

Whenever $M$ is a multiset of real numbers, we let $\prod M$ denote the
product of all the elements in $M$. Thus if $\la$ is a partition, then
$\prod\mathcal{H}(\la)$ is the product of all hook lengths in $\la$. 
Using the notation of Lemmas~\ref{betaparthook},~\ref{betasymb} we
have the following.

\begin{cor}
$\prod\mathcal{H}(X)=\prod\mathcal{H}(\la)=\prod\mathcal{H}^{\delta^*}(S)$. 
\end{cor}

\begin{cor}\label{partmetaprod} 
With the notation of Theorem~\ref{partmeta} we have
$$\prod \mathcal{H}(\la)=\prod\mathcal{H}(\la_{(d)})\cdot |\prod\mathcal{H}^{\delta}(Q)|.$$
\end{cor}

\medskip

The celebrated {\it hook formula} for the degrees of the irreducible
characters of the symmetric group $S_n$  may be formulated as follows.

\begin{theorem}\label{degree} 
Let $\la$ be a partition of $n$, and let
$\chi_{\la}$ be the irreducible character of $S_n$ labelled by $\la$. Then
$$\chi_{\la}(1)=\frac{n!}{\prod\mathcal{H}(\la)}\:.$$
\end{theorem}

\smallskip

Then Corollary~\ref{partmetaprod} may be formulated as follows:
\begin{cor} \label{managen} 
If $|\la_{(d)}|=r$ then with the notation of Theorem~\ref{partmeta} 
we have
$$
\chi_{\la}(1)=
\frac{n!}{r!}\frac{1}{|\prod\mathcal{H}^{\delta}(Q)|}\chi_{\la_{(d)}}(1)\:.$$
\end{cor}

\smallskip

\begin{rem} {\rm Corollory~\ref{managen} is equivalent to a generalization of
\cite[Theorem 9.1]{MaNa}. In this theorem, $d$ is assumed to be a
prime, the $\beta$-set $X$ for $\la$ is chosen to be the set of first column
  hook lengths for $\la$, and the short hooks in the balanced symbol
  $Q$ (which is called $S$ in \cite[Theorem 9.1]{MaNa}) are reversed. As pointed
  out in Remark~\ref{hookrev} above, the reversal of short hooks does not
  influence the absolute value of the products
  $\prod\mathcal{H}^{\delta}(Q)$
of all hook lengths.  }
\end{rem}

\bigskip
\section{Minimal data tuples and hooks in symbols} \label{symbs}

Let $S=(X_0,X_1,\ldots,X_{d-1}) $ be  a $d$-symbol.
Given $\ell \in \N$ and $e \in [d]$ we want to consider
the $(\ell,e)$-core and $(\ell,e)$-quotient of $S$. 
We start by the case $e=0$, where our result Theorem~\ref{ell0} is
slightly stronger than in the general case, Theorem~\ref{elle}. Only the
short hooks create a difficulty in the general case.

First we define $S_{*\ell}$ as the $d\ell$-symbol
$$S_{*\ell}:=s_{d\ell}(s_d^{-1}(S)).$$
Here $s_d^{-1}$ transforms the $d$-symbol $S$ into a $\beta$-set $X$
and $S_{*\ell}$ is then the $d\ell$-symbol associated to $X$. Thus
$S_{*\ell}$ may be seen as the ``splitting of $S$ into $\ell$ pieces''.

By Definition ~\ref{squotcore} the $d\ell$-symbol
$S_{*\ell}$ has a balanced quotient which we call
the {\it balanced $\ell$-quotient of $S:$}
$$Q_{\ell}(S):=Q(S_{*\ell}).$$
The  $d\ell$-symbol $S_{*\ell}$ also has  a core $C(S_{*\ell})$.
By Lemma~\ref{canbijs} there is a bijection between $H(S)$ and
$H(S_{*\ell})$  which may be described as follows. Let
$X=s_d^{-1}(S)$. Consider $z=(a,b,i,j) \in  H(S)$. Here $i,j \in
[d]$. Write $a=r\ell+s, b=r'\ell+s', ~s,s' \in [\ell]~$. 
Then $\mathfrak{h}_{X,d}^{-1}$ maps $(a,b,i,j)$ to
$(ad+i,bd+j)=(r(d\ell)+sd+i,r'(d\ell)+s'd+j) \in H(X)$ which by
$\mathfrak{h}_{X,d\ell}$ is mapped to $z'=(r,r',sd+i,s'd+j)$. 
This bijection restricts to a bijection between the
$(\ell,0)$-hooks in $H(S)$ and the $(1,0)$-hooks in
$H(S_{*\ell})$. More generally
it restricts to a bijection between the $(k\ell,0)$-hooks in $H(S)$
and the $(k,0)$-hooks in $H(S_{*\ell})$. 
This may also be applied to the $d$-symbol $C_{(\ell)}(S)$ having the
property that ${C_{(\ell)}(S)}_{*\ell}=C(S_{*\ell})$. 
(Thus $C_{(\ell)}(S)=s_d(s_{d\ell}^{-1}(C(S_{*\ell}))).$)
It shows that
$C_{(\ell)}(S)$ has no $(\ell,0)$-hooks and it is thus called the
{\it $\ell$-core of $S.$} It is really the $(\ell,0)$-core of $S$. 
The reader should notice that the quotient $Q_{\ell}(S)$ is a $d\ell$-symbol,
whereas the core $C_{(\ell)}(S)$ is a $d$-symbol. (Subscripts with
brackets are used for $d$-symbols and subscripts without
brackets for $d\ell$-symbols.)

\smallskip

Here is an example:
\begin{example} \label{ex3} \rm{ Choose $d=2, \ell=3$. 
Consider $S=(X_0,X_1)$ with $X_0=\{9,7,4,2\}$, $X_1=\{3,1,0\}$. 
We get $X=s_2^{-1}(S)=\{18,14,8,7,4,3,1\}$
and $p(S)=p(X)=(12,9,4,4,2,2,1)$.
\iffalse{=====
Putting $X$ on the $d\ell=6$-abacus we get
\begin{center}
\begin{tabular}{c c c c c c }
0      &\bf{1}&2      &\bf{3}&\bf{4}&5 \\
6      &\bf{7}&\bf{8} &9     &10    &11 \\
12     &13    &\bf{14}&15    &16    &17 \\
\bf{18}&19    & 20    &21    &22    &23
\end{tabular}
\end{center}
====}\fi
We compute $S_{*3}=(X'_0,X'_1,\ldots,X'_5)$, where
$$X'_0=\{3\},X'_1=\{1,0\},X'_2=\{2,1\},X'_3=\{0\},X'_4=\{0\},
X'_5=\emptyset $$
and get as balanced 3-quotient of $S$ the 6-symbol
$$Q_3(S)=t_6((3),(0),(1,1),(0),(0), (0))=(Y_0,Y_1,\ldots,Y_5)$$
where
$$Y_0=\{4,0\},Y_1=\{1,0\},Y_2=\{2,1\},Y_3=\{1,0\},Y_4=\{1,0\},
Y_5=\{1,0\}\:. $$
We have  $C(S_{*3})=([1],[2],[2],[1],[1],[0])$, 
and $s_6^{-1}(C(S_{*3}))=\{8,7,4,3,2,1,0\}$ so that
$C_{(3)}(S)$ is the 2-symbol
$$C_{(3)}(S)=(\{4,2,1,0\},\{3,1,0\}) \:.$$
}
\end{example}

\medskip
Let $\delta=(c_0,c_1,\ldots,c_{d-1};m)$ be a $d$-hook data
tuple. Consider the $d\ell$-hook data tuple
$\delta_{*\ell}=(c'_1,\ldots,c'_{d\ell-1};m \ell)$
where we have $c'_{i'}=sm+c_i$ if $i'=sd+i$. Thus
\begin{gather*}
\delta_{*\ell}=(c_0,c_1,\ldots,c_{d-1},\\
m+c_0,m+c_1,\ldots,m+c_{d-1},\\
2m+c_0,2m+c_1,\ldots, 2m+c_{d-1}\\
\ldots,\\
(\ell-1)m+c_0, \ldots, (\ell-1)m+c_{d-1};m\ell).
\end{gather*}
We then have that
if  $z\in H(S)$ is mapped to $z' \in H(S_{*\ell})$ then
$h^{\delta}(z)=h^{\delta_{*\ell}}(z')$. 
Indeed, if $z=(a,b,i,j)$ is mapped to  $z'=(r,r',sd+i,s'd+j)$ in the
above notation then
\begin{gather*}
h^{\delta}(z)=(a-b)m+c_i-c_j \\
=(r-r')m\ell+(sm+c_i)-(s'm+c_j)=h^{\delta_{*\ell}}(z').
\end{gather*}
Thus
$$\mathcal{H}^{\delta}(S)=\mathcal{H}^{\delta_{*\ell}}(S_{*\ell}).$$

This may be applied to $\delta=\delta^o=(0,0,\ldots, 0;1)$, the
minimal $d$-hook data tuple (Remark~\ref{choice}).
The length $h^{\delta^o}(z)$ of the hook $z=(a,b,i,j) \in H(S)$
is then $h(z)=a-b$. Thus
$\mathcal{H}(S):=\mathcal{H}^{\delta^o}(S)$ is just the
multiset of hook lengths of $S$, as defined in \cite{Ma}, including short
hooks of length 0.

We have
$$\delta^o_{*\ell}=(0,\ldots,0,1,\ldots,1,\ldots,\ell-1,\ldots,\ell-1;
\ell)\:.$$
Now apply
Theorem~\ref{metasatz} to the $d\ell$-symbol $S_{*\ell}$ and $\delta^o_{*\ell}$. 
The balanced quotient $Q$ in the theorem is just the balanced
$\ell$-quotient of $S_{*\ell}$, i.e., 
$Q=Q_{\ell}(S)$. The core in the theorem is $C(S_{*\ell})$. By
definition of $C_{(\ell)}(S)$ we have
$C_{(\ell)}(S)_{*\ell}=C(S_{*\ell})$. 
Corresponding to $\delta_S$ in the theorem, we have here  
\begin{gather*}
\delta_{\ell,S}=(x_{0,0}\ell,\ldots,x_{d-1,0}\ell,\\
1+x_{0,1}\ell,\ldots,1+x_{d-1,1}\ell,\\
\ldots,\\
(\ell-1)+x_{0,\ell-1}\ell,\ldots,(\ell-1)+x_{d-1,\ell-1}\ell;\ell),
\end{gather*}
where $x_{i,j}$ is the number of elements in
$X_i$ which are congruent to $j$ modulo~$\ell$. 
By the above
$\mathcal{H}^{\delta^o}(S)=\mathcal{H}^{\delta^o_{*\ell}}(S_{*\ell})$
and
$\mathcal{H}^{\delta^o}(C_{(\ell)}(S))=\mathcal{H}^{\delta^o_{*\ell}}(C(S_{*\ell}))$. 
As in Theorem~\ref{betameta} we have that the hook
lengths are non-negative. 
Thus we obtain 

\begin{theorem}\label{ell0} 
Let $S=(X_0,X_1,\ldots,X_{d-1})$ be a $d$-symbol.
  %and that for $i \in [d], j\in [r]$
Let $C=C_{(\ell)}(S)$ be the $\ell$-core of $S$ and let $Q=Q_{\ell}(S)$ 
be the balanced $\ell$-quotient of $S$.  
Then
$$\mathcal{H}(S)=abs(\mathcal{H}^{\delta_{\ell,S}}(Q)) \cup
\mathcal{H}(C)$$
where $abs(\mathcal{H}^{\delta_{\ell,S}}(Q))$ is the multiset of all
$|h^{\delta_{\ell,S}}(z)|$, $z \in H(Q)$.
\end{theorem}

We continue Example~\ref{ex3}.
\begin{example} \label{ex4}{\rm 
We have $d=2,\ell=3$,  $S=(X_0,X_1)$ with $X_0=\{9,7,4,2\}$, $X_1=\{3,1,0\} $ 
and $p(S)=(12,9,4,4,2,2,1)$. 
The hook lengths (including 3 short hooks of length 0) of $S$ may
conveniently be read off the
2-abacus for $X=s_2^{-1}(S)$ and recorded systematically in the
Young diagram of $p(S)$.
The 2-abacus representation of $X$ 
(where the subscripts indicate the elements of
$X_0$ and  $X_1$) is 
\begin{center}
\begin{tabular}{c c} 
0&\bf{1}$_0$ \\
2&\bf{3}$_1$ \\
\bf{4}$_2$&5 \\
6&\bf{7}$_3$\\
\bf{8}$_4$&9\\
10&11\\
12&13\\
\bf{14}$_7$&15\\
16&17\\
\bf{18}$_9$&19
\end{tabular}
\end{center}
The hook lengths of $S$, i.e., the elements of $\mathcal{H}(S):$
\begin{center}
\begin{tabular}{c c c c c c c c c c c c}
9&8&7&6&5&4&4&3&3&2&1&1\\
7&6&5&4&3&2&2&1&1&&&\\
4&3&2&1\\
3&2&1&0\\
2&1\\
1&0\\
0
\end{tabular}
\end{center}
We calculated $C=C_{(3)}(S)=(\{4,2,1,0\},\{3,1,0\}) $ and get
$p(C)=(2,2)$ with hook lengths
\begin{center}
\begin{tabular}{c c }
2&1\\
1&0
\end{tabular}
\end{center}
so that $\mathcal{H}(C)=\{2,1,1,0\}$.
We also calculated $Q=Q_3(S)=(Y_0,Y_1,\ldots,Y_5)$ where
$Y_0=\{4,0\},Y_1=\{1,0\},Y_2=\{2,1\},Y_3=\{1,0\},Y_4=\{1,0\},
Y_5=\{1,0\}. $ We put $Q$ on the 6-abacus

\begin{center}
\begin{tabular}{c c c c c c }
\bf{0}$_0$ &\bf{1}$_0$&2      &\bf{3}$_0$&\bf{4}$_0$ &\bf{5}$_0$ \\
6      &\bf{7}$_1$&\bf{8}$_1$ &\bf{9}$_1$&\bf{10}$_1$&\bf{11}$_1$ \\
12     &13    &\bf{14}$_2$&15    &16    &17 \\
18     &19    &20     &21    &22    &23 \\
\bf{24}$_4$&25    & $\cdots$
\end{tabular}
\end{center}
The hook data tuple $\delta_{3,S}$ in Theorem~\ref{ell0} is here
$\delta_{3,S}=(3,6,7,4,5,2;3)$. 
We record the $\delta_{3,S}$-hook lengths of $Q$ in the Young diagram of

\noindent $p(Q)=(13,4,2,2,2,2,2,1,1,1)$
\begin{center}
\begin{tabular}{ c c c c c c c c c c c c c}

 8&9&6&3&5&4&7&3&0&-1&2&1&4\\
 6&7&4&1&\\
-2&-1\\
 1&2\\
 0&1\\
 3&4\\
 2&3\\
-5\\
-2\\
-3
\end{tabular}
\end{center}
%\smallskip
The set of absolute values of these elements form
$abs(\mathcal{H}^{\delta_{3,S}}(Q))$ and we see that 
$\mathcal{H}(S)=abs(\mathcal{H}^{\delta_{3,s}}(Q)) \cup
\mathcal{H}(C)$. }
  \end{example}
%\medskip
We want to generalize this to $(\ell,e)$-cores and
$(\ell,e)$-quotients, $e \in [d]$, the above Theorem~\ref{ell0}
being the case $e=0$. 
Involving a suitable permutation of $\N_0$, 
this is done easily. 
This permutation $\sigma=\sigma_{d,\ell,e}$ is defined as follows.
Any $n \in \N_0$ may be written uniquely as
$n=r(d\ell)+sd+t,~~r\in \N_0,~s \in [\ell], ~t \in [d]$. 
Then $\sigma_{d,\ell,e}(n)=r(d\ell)+sd+(t+re)_{[d]}$.
The $(\ell,e)$-{\it twist} $\sigma(S)$ of the $d$-symbol $S$ is defined as
$
\sigma(S)=s_{d}\sigma(s_d^{-1}(S)).
$
Note that if $S=(X_0,X_1,\ldots,X_{d-1})$, then the $d$-symbol $\sigma(S)$ is
{\it not} equal to $(\sigma(X_0),\sigma(X_1), \ldots, \sigma(X_{d-1}))$. But if $X=s_d^{-1}(S)$
is the $\beta$-set of $S$, then $\sigma(X)$ is the $\beta$-set of $\sigma(S)$.

Let us consider a {\it long} hook $z=(a,b,i,j) \in H(S)$, i.e., we
have $a>b$. Write $a=r\ell+s, b=r'\ell+s', ~s,s' \in [\ell]~$. Then
\begin{center}
$\xymatrix{(a,b,i,j) \in H(S) \ar[d]^-{\mathfrak{h}_{X,d}^{-1}} \\
  (ad+i,bd+j)=(r(d\ell)+sd+i,r'(d\ell)+s'd+j) \in H(X)
  \ar[d]^-{\sigma} \\
(\sigma(ad+i),\sigma(bd+j))=}$
$\xymatrix{ (r(d\ell)+sd+(i+re)_{[d]},r'(d\ell)+s'd+(j+r'e)_{[d]}) \in H(\sigma(X)) \ar[d]^-{\mathfrak{h}_{X,d}} \\
z'=(a,b,(i+re)_{[d]}, (j+r'e)_{[d]}) \in H(\sigma(S))
 } $
\end{center}
Here it should be noted that since $z$ is long, then
$\sigma(ad+i)>\sigma(bd+j)$ so that $(\sigma(ad+i),\sigma(bd+j)) \in
H(\sigma(X))$. Therefore the above is a bijection between the sets
 $H_{>0}(S)$ and $H_{>0}(\sigma(S))$ of {\it long} hooks.
There is in general no bijection for short hooks; their number may
differ!

It follows from the definitions that if
$\delta^o=(0,0,\ldots, 0;1)$ is the
minimal $d$-hook data tuple, then corresponding {\it long} hooks
in $S$ and $\sigma(S)$
have the same $\delta^o$-length. Now $\delta^o$-lengths are always
non-negative and the long hooks are exactly those of non-zero
$\delta^o$-length.
Thus
\begin{center}
$\mathcal{H}_{>0}^{\delta^o}(S)=\mathcal{H}_{>0}^{\delta^o}(\sigma(S))$. \end{center}

If  $z=(a,b,i,j) \in H(S)$ is an $(\ell,e)$-hook,
i.e., $a-b=\ell$, $(j-i)_{[d]}=e$,  then in the above notation $r'=r-1, s=s'$ and
$(a,b,i,j)$ is mapped to $z'=(a,b,(i+re)_{[d]},(j+(r-1)e)_{[d]}) \in H(\sigma(S)$.
Since
$$(j+(r-1)e)_{[d]}-(i+re)_{[d]})_{[d]}=(j-i-e)_{[d]}=0_{[d]}$$ we see that $z'$ is an
$(\ell,0)$-hook in  $H(\sigma(S))$ .

This shows that the $d$-symbol $C_{(\ell,e)}(S)$ satisfying that
$\sigma(C_{(\ell,e)}(S))=C_{(\ell)}(\sigma(S))$
is obtained from $S$ by removing all $(\ell,e)$-hooks.
(Thus $C_{(\ell,e)}(S)=s_d\sigma^{-1}s_d^{-1}(C_{(\ell)}(\sigma(S)))=
s_d\sigma^{-1}s_{d\ell}^{-1}(C(\sigma(S)_{*\ell})). $)
We call
 $C=C_{(\ell,e)}(S)$ the $(\ell,e)$-{\it core} of $S$.
We also need to define a (balanced) $(\ell,e)$-quotient
$Q=Q_{\ell,e}(S)$ of $S$. To do this we consider
the $d\ell$-symbol $S_{*\ell,e}$
defined by
\begin{center}
$S_{*\ell,e}:=s_{d\ell}\sigma(s_d^{-1}(S))$.
\end{center}
Thus $S_{*\ell,e}$ is the  $d\ell$-symbol associated to the
$\beta$-set $\sigma(X)$.
If again $\sigma(S)=s_{d}\sigma(s_d^{-1}(S))$ is the $(\ell,e)$-twist of $S$, then by definition
\begin{center}
$\sigma(S)_{*\ell}=S_{*\ell,e}$.
\end{center}
We define
 $Q=Q_{\ell,e}(S):=Q(S_{*\ell,e})=Q(\sigma(S)_{*\ell})$, which is a $d\ell$-symbol.

We now apply Theorem~\ref{ell0} to the $d$-symbol $\sigma(S)$. 
We get for a suitable $d\ell$-hook data tuple
$\delta=\delta_{\ell,\sigma(S)}$, defined in analogy with $\delta_{\ell,S}$ above,
$$\mathcal{H}(\sigma(S))=abs(\mathcal{H}^{\delta}(Q)) \cup
\mathcal{H}(C),$$
where $abs(\mathcal{H}^{\delta}(Q))$ is the multiset of all
$|h^{\delta}(z)|$, $z \in H(Q)$.
In this multiset equality we remove all occurrences of 0 and get
$$\mathcal{H}_{>0}(\sigma(S))=abs(\mathcal{H}_{>0}^{\delta}(Q)) \cup
\mathcal{H}_{>0}(\sigma(S)_{(\ell)}).$$
By the above $\mathcal{H}_{>0}(\sigma(S))=\mathcal{H}_{>0}(S)$ and
$\mathcal{H}_{>0}(\sigma(S)_{(\ell)})=\mathcal{H}_{>0}(S_{(\ell,e)})$. We have shown

\begin{theorem}\label{elle} Suppose that $S=(X_0,X_1,\ldots,X_{d-1})$
  is a $d$-symbol. Let $C=C_{(\ell,e)}(S)$ be
  the $(\ell,e)$-core of $S$, $Q=Q_{\ell,e}(S)$ the balanced
  $(\ell,e)$-quotient of $S$ and $\delta=\delta_{\ell,\sigma(S)}$.  Then
$$\mathcal{H}_{>0}(S)= abs(\mathcal{H}_{>0}^{\delta}(Q)) \cup
\mathcal{H}_{>0}(C),$$
where $abs(\mathcal{H}_{>0}^{\delta}(Q))$ is the multiset
of all non-zero
$|h^{\delta}(z)|$, $z \in H(Q)$.
\end{theorem}

\begin{rem} {\rm In analogy with Corollary \ref{managen} we may use the
    Theorems \ref{ell0} and \ref{elle} to prove relative hook formulas
    for unipotent degrees. See \cite[Proposition (8.9)]{OlL} for
    unipotent character degrees in finite classical groups or more
    generally \cite[3.12]{Ma}. }
\end{rem}

\begin{example} \label{ex5}{\rm Let us consider the case  $d=2,\ell=3,
    e=1$. The permutation $\sigma=\sigma_{2,3,1}$ is a product of
    transpositions:
$$\sigma=\Pi_{k ~odd}(6k,6k+1)(6k+2,6k+3)(6k+4,6k+5).$$
We want to apply Theorem~\ref{elle} to the 2-symbol $S=(X_0,X_1)$ 
with $X_0=\{9,7,4,2\}$, $X_1=\{3,1,0\}$.  
We have $X=s_2^{-1}(S)=\{18,14,8,7,4,3,1\}$. 
Thus $\sigma(X)=\{19,14,9,6,4,3,1\}$ and
$$\sigma(S)=s_2(\sigma(X))=(\{7,3,2\},\{9,4,1,0\})\:.$$
Putting $\sigma(S)$ on the
$d\ell=6$-abacus we get in analogy with Example~\ref{ex3} that
$$\sigma(S)_{*3}=S_{*3,1}=(\{1\},\{3,0\},\{2\},\{1,0\},\{0\},\emptyset)$$
and the balanced quotient is
$$Q=Q_{3,1}(S)=(\{1\},\{2\},\{2\},\{0\},\{0\},\{0\}).$$
We have $C(\sigma(S)_{*3})=([1],[2],[1],[2][1],\emptyset)$. Thus $\sigma(S)$ has 3-core
$$C=C_{(3,1)}(S)=s_2\sigma^{-1}s_6^{-1}(C(\sigma(S)_{*3}))=([5],[2])$$ with
  $$\mathcal{H}_{>0}(C)=\{2,1,1\}.$$
We have $\delta_{\ell,\sigma(S)}=(3,6,4,7,5,2;3)$ and
list the $\delta_{\ell,\sigma(S)}$-hooks of $Q$ in the Young
diagram of $p(Q)=(9,9,3,3,3,3):$
\begin{center}
\begin{tabular}{ c c c c c c c c c }

 7&4&6&1&3&0&2&5&1\\
 9&6&8&3&5&2&4&7&3\\
 3&0&2\\
-1&-4&-2\\
 2&-1&1\\
 4&1&3

\end{tabular}
\end{center}

The union of the multiset of the non-zero absolute values of these hook lengths
and  $\mathcal{H}_{>0}(C)$ coincides with the multiset of hook lengths
of $S$ listed in the Young diagram of
$p(S)$ in Example~\ref{ex4}. 
This is in accordance with Theorem~\ref{elle}.
}
\end{example}

\noindent\textbf{Acknowledgement.}\quad
The authors wish to thank Gabriel Navarro for a 
discussion leading to this investigation.

\end{document}